\documentclass[a4paper,12pt]{article}
\usepackage{amsmath}
\usepackage{amsthm}
\usepackage{graphicx}
\usepackage{verbatim}
\usepackage{epsfig}
\usepackage{hyperref}
\usepackage{float}
\usepackage{color}
\usepackage{fullpage}
\usepackage{amsfonts}
\usepackage{amscd}
\usepackage{amssymb}
\usepackage{graphics}
\usepackage{amsmath}
\usepackage[english]{babel}
\usepackage{bbm}
\usepackage{yfonts}
\usepackage{mathrsfs}
\usepackage{color}
\DeclareFontFamily{OML}{rsfs}{\skewchar\font'177}
\DeclareFontShape{OML}{rsfs}{m}{n}{ <5> <6> rsfs5 <7> <8> <9>
rsfs7 <10> <10.95> <12> <14.4> <17.28> <20.74> <24.88> rsfs10 }{}
\DeclareMathAlphabet{\mathfs}{OML}{rsfs}{m}{n}

\newcommand{\BE}{{\mathbb{E}}}

\newcommand{\BP}{{\mathbb{P}}}

\newcommand{\BR}{{\mathbb{R}}}

\newcommand{\BZ}{{\mathbb{Z}}}

\newcommand{\CC}{{\mathcal{C}}}

\newcommand{\CM}{{\mathcal{M}}}

\newcommand{\CX}{{\mathcal{X}}}

\newtheorem{theorem}{Theorem}[section]
\newtheorem{proposition}[theorem]{Proposition}
\newtheorem{lemma}[theorem]{Lemma}

\begin{document}
\numberwithin{equation}{section} \numberwithin{figure}{section}
\title{The existence phase transition for two Poisson random fractal models.}
\author{Erik I. Broman\footnote{Department of Mathematics, Chalmers University of Technology and Gothenburg University, Sweden. E-mail: broman@chalmers.se. Research supported by the Swedish research Council}
 \ Johan Jonasson\footnote{Department of Mathematics, Chalmers University of Technology and Gothenburg University, Sweden. E-mail: jonasson@chalmers.se. Research supported by the Knut and Alice Wallenberg foundation}
 \ Johan Tykesson\footnote{Department of Mathematics, Chalmers University of Technology and Gothenburg University, Sweden. E-mail: johant@chalmers.se. Research supported by the Knut and Alice Wallenberg foundation}}
\maketitle
\begin{abstract}
In this paper we study the existence phase transition of the
random fractal ball model and the random fractal box model.
We show that both of these are in the empty phase at the critical
point of this phase transition.
\end{abstract}
\section{Introduction}

In order to better explain the rest of the paper,
we shall start by a rather informal description of the general setup
(see for example \cite{BC} for details).
Let $M$ be the set of bounded subsets of $\BR^d$ with
non-empty interior, and let $\CM$ be some (suitable) $\sigma$-algebra
on $M.$ We consider a measure $\mu$ on $(M,\CM)$ which is scale invariant
in the following sense. If $A\in \CM$ is such that $\mu(A)<\infty,$ then
$\mu(A_s)=\mu(A)$ where $0<s<\infty$ and
\[
A_s:=\{K:K/s\in A\}.
\]
We will also assume that $\mu$ is translation invariant in that 
$\mu(x+A)=\mu(A)$ for every $A \in \CM.$ Here of course, 
$x+A=\{L\subset \BR^d:L=x+K \textrm{ for some } K\in A\}$.

In order to define a model which will exhibit a non-trivial behaviour,
it is often necessary to restrict $\mu$ to sets of diameter smaller than
some cutoff. This is indeed what we do in this paper (see also the discussion 
in Section \ref{sec:defmodels}). For such measures, the property 
$\mu(A_s)=\mu(A)$ will still hold, but only if neither $A$ nor 
$A_s$ contains sets with diameter larger than the cutoff. We shall call 
such a measure {\em semi} scale invariant.
In the rest of this introduction, any measure $\mu$ we refer to will be
semi scale invariant.

Using $\lambda \mu$ where $0<\lambda<\infty$ as the intensity measure,
one can define a Poisson process $\Phi_\lambda(\mu)$ on $M$.
Thus constructed, $\Phi_\lambda(\mu)$ is a semi scale and translation 
invariant random collection
of bounded sets of $\BR^d.$ This setup contains many interesting examples such
as the Brownian loop soup introduced in \cite{LW}, and the semi scale invariant
Poisson Boolean model studied for instance in \cite{BE} (see also the
references therein). Throughout, this latter model will be referred to
simply as the {\em fractal ball model}, and we shall
give an exact definition of it in Section \ref{sec:defmodels}.
In this fractal ball model, the measure $\mu$ above is supported on the
set of {\em open}
balls of $\BR^d$. Of course, one could also consider a process of closed balls,
or indeed a mix of open and closed balls. As we will see, the results of
this paper
are also valid for these cases, see further the remark after the statement
of Theorem \ref{thm:main}.

Throughout this paper, we will let
\begin{equation} \label{eqn:defCC}
\CC(\Phi_\lambda(\mu))
:=\BR^d \setminus \bigcup_{K\in \Phi_\lambda(\mu)} K,
\end{equation}
and we will usually write $\CC(\lambda)$ or simply $\CC.$ Thus, with $\mu$
as above, $\CC$ is a semi scale invariant random fractal and we will be concerned
by various properties of $\CC(\lambda)$ as $\lambda$ varies. It is useful to observe that 
by using a standard coupling, $\CC(\lambda)$ is decreasing in $\lambda.$

Random fractal models exhibits several phase transitions (see for instance \cite{DM}).
However, the perhaps two most natural are the {\em existence} and the
{\em connectivity} phase transitions as we now explain.
Define
\[
\lambda_e:=\inf\{\lambda>0:\BP(\CC(\lambda)=\emptyset)=1\}.
\]
Therefore, for $\lambda>\lambda_e,$ $\CC(\lambda)$ is almost surely empty,
and we say that it is in the empty phase. If instead $\lambda<\lambda_e,$ then
$\BP(\CC(\lambda)\neq \emptyset)=1$. We say that $\lambda_e$ is the critical
point of the existence phase transition.
Analogously, we can define
\[
\lambda_c:=\sup\{\lambda>0:\BP(\CC(\lambda)
\textrm{ contains connected components larger than one point})=1\}.
\]
Thus, for $\lambda>\lambda_c,$ $\CC(\lambda)$ is almost surely
totally disconnected, while for $\lambda<\lambda_c,$ $\CC(\lambda)$
will contain connected components.

Of course, whenever such phase transitions occur, it is natural and
interesting to ask what happens at the critical points. In \cite{BC}
it was proven in full generality that
\[
\BP(\CC(\lambda_c)
\textrm{ contains connected components larger than one point})=1,
\]
so that {\em at} $\lambda_c$ the fractal is in the connected phase.
Thus, this phase transition is very well understood.

The existence phase transition is much less understood.
Hitherto, the only exact results appear to be in dimension 1. Indeed,
in \cite{Shepp}, exact conditions for when random intervals
cover a line were established. However, there has been some
progress (see \cite{BE}) on the case of the fractal ball model in $d\geq 2$,
see Section \ref{sec:defmodels} for a precise statement of these results.

In analogy with how the fractal ball model is defined, we can also define the
{\em fractal box model} (again see Section \ref{sec:defmodels}) for
which the measure $\mu$ is supported on boxes of the form $(a,b)^d$
for $a<b.$ In this case, $\Phi_\lambda(\mu)$ is then a random semi scale
invariant collection of boxes in $\BR^d.$
Whenever we need to distinguish between the ball and the box model,
we shall write $\CC^{ball}$ and $\CC^{box}$ etc.

Let $v_d$ be the volume of the unit ball in $\BR^d.$
The main result of this paper is the following.
\begin{theorem} \label{thm:main}
For any $d\geq 1,$ we have that $\lambda_e^{ball}=d/v_d$ while
$\lambda^{box}_e=d.$ Furthermore,
\[
\BP(\CC^{box}(\lambda^{box}_e)=\emptyset)
=\BP(\CC^{ball}(\lambda^{ball}_e)=\emptyset)=1.
\]
\end{theorem}
\noindent
{\bf Remarks:} The fact that $\lambda_e^{ball}=d/v_d$ is easily deduced
from results in \cite{BE}, while we determine $\lambda_e^{box}$
by a straightforward second moment argument.
Thus, the main contribution of this paper is to determine what happens
{\em at} the critical point of these phase transitions.

If we choose to consider closed balls (boxes) in place of open, then of course
we would have that $\CC_{closed}\subset \CC_{open}$ (using obvious notation).
However, when determining $\lambda_e$, one sees that the argument does not
depend on whether we use open or closed sets so that
$\lambda_e(\CC_{closed})=\lambda_e(\CC_{open}).$ It then follows trivially
that Theorem \ref{thm:main} holds also for the case of closed balls (boxes).

The result does not depend on the specific value of the cutoff (as is clear
from the proofs). However, it requires some cutoff. 
\bigskip

The rest of the paper is organized as follows. In Section \ref{sec:defmodels}
we give precise definitions of our models and also provide some further
background. In Section \ref{sec:proofs}, we will prove Theorem \ref{thm:main}.

\section{Models} \label{sec:defmodels}

We start by defining the fractal ball model, although we will later reuse 
much of the notation for the box model.

Let $\nu$ be a locally finite measure on $(0,1]$,
and let $\mu=dx \times \nu$ (where $dx$ denotes $d$-dimensional
Lebesgue measure) denote the resulting product measure on $\BR^d\times (0,1].$
Then, we let $\Phi_\lambda(\mu)$ be a Poisson process on $\BR^d\times (0,1]$
using $\lambda \mu$ as the intensity measure. This definition might seem
to clash with $\Phi_\lambda(\mu)$ defined in the introduction (which
was a Poisson process on sets). However, this is easily resolved
by associating the point $(x,r)\in\BR^d\times (0,1]$ with the open
ball $B(x,r)$ centered at $x$ and with radius $r.$ Thus, we might write
\eqref{eqn:defCC} as
\[
\CC(\Phi_\lambda(\mu))
:=\BR^d \setminus \bigcup_{(x,r)\in \Phi_\lambda(\mu)} B(x,r).
\]
Let $A:=\{(x,r)\in\BR^d\times [\epsilon,1]: o\in B(x,r)\}$, and let
$A_{\epsilon^{-1}}=\{(x,r)\in \BR^d\times [\epsilon^2,\epsilon]: o\in B(x,r)\}$
(where $o$ denotes the origin).
We observe that if $\nu(dr)=r^{-d-1}dr,$ then we have that (with $I(\cdot)$
being an indicator function)
\[
\mu(A)=\int_\epsilon^1 \int_{\BR^d} I(|x|\leq r) dx\nu(dr)
=v_d \int_\epsilon^1 r^d r^{-d-1}dr=-v_d\log \epsilon,
\]
and an analogous calculation shows that also
$\mu(A_{\epsilon^{-1}})=-v_d\log \epsilon$. 

We observe that $\mu$ cannot be fully scale invariant since  
we have that $\mu(A_\epsilon)=0.$ This follows since 
$A_\epsilon$ only contains sets with $r\geq 1.$
However, if we in the above replace $\nu$ by  $\tilde{\nu}(dr)=r^{-d-1}dr$ 
supported on $(0,\infty)$, we would obtain a fully scale invariant measure 
$\tilde{\mu}.$ Thus, our measure $\mu$ is the restriction of $\tilde{\mu}$
to sets with $r\leq 1,$ which is then our cutoff. 
In particular we have that $\mu(A)=\mu(A_s)$
as long as neither $A$ nor $A_s$ contains sets with $r>1.$ 
We note that it would perhaps be more proper to write 
$\nu(dr)=I(0<r\leq 1)r^{-d-1}dr$. However, we will allow ourselves to 
slightly abuse notation by writing $\nu(dr)=r^{-d-1}dr,$ and remembering that 
$\nu$ is supported on $(0,1]$.

It is certainly possible to consider other choices of $\nu$, but in this paper
we shall focus on the semi scale invariant case. 
However, we want to mention the following 
result from \cite{BE} which deals with other choices of $\nu.$

\begin{theorem}[From \cite{BE}] \label{thm:BE} For the fractal ball model,
if $\BP(\CC=\emptyset)=1$ then
\begin{equation} \label{eqn:nec}
\int_0^1 u^{d-1}
\exp\left(\lambda v_d \int_u^1 r^{d-1}(r-u)\nu(dr)\right)du=\infty,
\end{equation}
while if
\begin{equation} \label{eqn:suff}
\limsup_{u\to 0}u^d\exp\left(\lambda v_d\int_u^1 (r-u)^d\nu(dr)\right)du=\infty,
\end{equation}
then $\BP(\CC=\emptyset)=1$.
\end{theorem}
{\bf Remark:} 
Taking $\nu(dr)=r^{-d-1}dr,$ one concludes from
\eqref{eqn:nec} and \eqref{eqn:suff} that $\lambda^{ball}_e=d/v_d$.
However, simple calculations reveal that \eqref{eqn:suff} is not
satisfied for $\lambda=d/v_d,$ and so
we cannot conclude whether $\BP(\CC(\lambda^{ball}_e)=\emptyset)=1$.
As pointed out in \cite{BE}, it follows from Theorem \ref{thm:BE}
that if $\nu(dr)=r^ {-d-1}(1+2|\log (r)|^{-1})$ and $\lambda=d/v_d,$ then
$\BP(\CC=\emptyset)=1$ while if
$\nu(dr)=r^ {-d-1}(1-2|\log (r)|^{-1})$ and $\lambda=d/v_d,$ then
$\BP(\CC=\emptyset)=0$. Thus, although their results do not cover
the critical case, it comes logarithmically close. Of course,
Theorem \ref{thm:main} improves on Theorem \ref{thm:BE}
in that we here determine the critical case.

\bigskip

We now turn to the fractal box model. Here, we again use the measures
$\nu$ and $\mu$ as above, but to any $(x,r)\in \BR \times (0,1],$
we associate the open box $X(x,r):=x+(-r/2,r/2)^d$. We then write
\[
\CC(\Phi_\lambda(\mu))
:=\BR^d \setminus \bigcup_{(x,r)\in \Phi_\lambda(\mu)} X(x,r).
\]
Letting $A=\{(x,r)\in \BR \times [\epsilon,1]:
X(x,r)\cap [0,1]^d\neq \emptyset\}$
we have that
\[
\mu(A)=\int_{\epsilon}^1\int_{\BR^d} I(x\in (-r/2,1+r/2)^d) dx\nu(dr)
=\int_{\epsilon}^1(1+r)^dr^{-d-1}dr.
\]
Similarly, if $A_{\epsilon^{-1}}
=\{(x,r)\in \BR \times [\epsilon^2,\epsilon]:
X(x,r)\cap [0,\epsilon]^d\neq \emptyset\}$ then
\[
\mu(A_{\epsilon^{-1}})
=\int_{\epsilon^2}^\epsilon\int_{\BR^d} I(x\in (-r/2,\epsilon +r/2)^d)dx \nu(dr)
=\int_{\epsilon^2}^\epsilon (\epsilon+r)^d r^{-d-1}dr
=\mu(A),
\]
so that also this model is semi scale invariant.

Whenever convenient, we will write $K\in \Phi$ to mean either a ball or a box,
depending on the context.

\section{Proofs} \label{sec:proofs}
We start this section by introducing some useful notation. First, let
\[
\bar{\CX}_n:=\left\{x+[0,1/n]^d:
x\in \left(\frac{1}{n}\BZ^d\right)\cap[0,1-1/n]^d\right\}.
\]
If $\bar{X}\in \bar{\CX}_n,$ we shall refer to $\bar{X}$ as a level $n$ box.
Note that the members $\bar{X}$ of $\bar{\CX}_n$ are deterministic, closed boxes. 
These should not be confused with the open boxes $X(x,r)$ that belong to the 
Poisson process $\Phi_\lambda^{box}$.

The interpretation of the following definitions differ
depending on whether we are considering the ball model or the box model.
However, we believe that this should not lead to any confusion.
For these models, we let
\[
\Phi_n:=\{(x,r)\in \Phi_\lambda (\mu):1/n\leq r \leq 1\},
\]
and define
\[
\CC_n:=\BR^d \setminus \bigcup_{K \in \Phi_{n}} K
\]
(where $K$ is then either a ball or a box). Thus, $\CC_n \downarrow \CC.$
For $m>n,$ let
\[
\CC_{m}^n:=\BR^d \setminus \bigcup_{K \in \Phi_{m}\setminus \Phi_n} K,
\]
so that $\CC_m^n \cap \CC_n=\CC_m$, and $\CC_m^n, \CC_n$ 
are independent. 
For any integer $n,$ let
\[
M_n:=\{\bar{X}\in \bar{\CX}_n:
\not \exists K \in \Phi_n\,:\,\bar{X} \subset K \}.
\]
Thus, $M_n$ is the set of level $n$ boxes which are not covered by a single
set in the Poisson process $\Phi_n.$ Then, let
\[
m_n:=\{\bar{X}\in \bar{\CX}_n:
\not \exists K \in \Phi_n\,:\,\bar{X} \cap K\neq \emptyset \},
\]
which is the set of level $n$ boxes untouched by the Poisson process $\Phi_n$.
We see that if $\bar{X}\in m_n,$ then in fact $\bar{X} \subset \CC_n$. 
Obviously, $|m_n|\leq |M_n|$ since an untouched box cannot be covered.

The following proposition is a part of Theorem \ref{thm:main}.
\begin{proposition} \label{prop:crit}
For the box model we have that $\lambda_e\geq d$.
\end{proposition}
\noindent
{\bf Proof.}
We start by noting that if $m_n\neq \emptyset$ for infinitely
many $n\geq 1,$  then
$\CC_n\cap[0,1]^d\neq \emptyset$ for every $n\geq 1.$ Since
$\CC_n\supset \CC_{n+1}$ for every $n,$ and the sets
$\CC_n\cap[0,1]^d$ are compact, we must then have that
\[
\CC\cap[0,1]^d =\bigcap_{n=1}^\infty\CC_n\cap[0,1]^d\neq \emptyset.
\]
We will prove that for $\lambda<d,$ there exists $c=c(\lambda)>0$ such that
\begin{equation} \label{eqn:mnpos}
\BP(m_n>0)\geq c,
\end{equation}
for every $n\geq 1.$ Then, we can conclude that
\[
\BP(m_n>0 \textrm{ infinitely often})
\geq \limsup_{n \to \infty}\BP(m_n>0)\geq c,
\]
by the reverse Fatou's lemma.

We shall proceed by proving \eqref{eqn:mnpos} using a second moment argument.
To that end, observe that by translation invariance, for any $\bar{X}\in \bar{\CX}_n,$ 
\begin{eqnarray} \label{eqn:ref1}
\lefteqn{\BP(\bar{X}\in m_n)=\BP([0,1/n]^d\in m_n)}\\
& & =\exp\left(-\lambda \mu(\{(x,r)\in \BR^d\times[1/n,1]:\ [0,1/n]^d \cap X(x,r)\neq \emptyset\})\right) \nonumber \\
& & =\exp\left(-\lambda\int_{1/n}^1 \int_{\BR^d} I(x\in(-r/2,r/2+1/n)^d)
dx \nu(dr)\right)\nonumber \\
& & =\exp\left(-\lambda\int_{1/n}^1 (r+1/n)^d r^{-d-1}dr\right) \nonumber \\
& & =\exp\left(-\lambda\int_{1/n}^1 r^{-d-1}
\sum_{k=0}^d {d \choose k} r^k n^{k-d} dr\right)\nonumber \\
& & =\exp\left(-\lambda\sum_{k=0}^d {d \choose k} n^{k-d}
\int_{1/n}^1 r^{k-d-1}dr\right)\nonumber \\
& & =\exp\left(-\lambda\log n-
\lambda\sum_{k=0}^{d-1} {d \choose k} n^{k-d}
\left(\frac{n^{d-k}-1}{d-k}\right)\right).\nonumber 
\end{eqnarray}
Since
\begin{equation}\label{eqn:ref2}
0\leq \sum_{k=0}^{d-1} {d \choose k} n^{k-d}\left(\frac{n^{d-k}-1}{d-k} \right)
\leq \sum_{k=0}^{d}{d \choose k}=2^d,
\end{equation}
we conclude that
\begin{equation} \label{eqn:Pmn}
e^{-\lambda 2^d}n^{-\lambda}\leq \BP(\bar{X}\in m_n)\leq n^{-\lambda}.
\end{equation}
Therefore,
\begin{equation} \label{eqn:fmoment}
\BE[m_n]=n^d\BP([0,1/n]^d\in m_n)\geq e^{-\lambda 2^d} n^{d-\lambda}.
\end{equation}

For $\bar{X}_1,\bar{X}_2\in \bar{\CX}_n$ let
$R^n_i:=\{(x,r)\in \BR^d \times [1/n,1]:\bar{X}_i \cap X(x,r)\neq \emptyset\}$
for $i=1,2.$ We have that
$\mu(R_1^n \cup R_2^n)=2\mu(R_1^n)-\mu(R_1^n \cap R_2^n)$.
First, we observe that 
\begin{eqnarray} \label{eqn:muR1}
\lefteqn{\mu(R_1^n) =\int_{1/n}^1 \int_{\BR^d} I(X(x,r)\cap \bar{X}_1 \neq \emptyset)  
dx \nu(dr) }\\
& & =\int_{1/n}^1 (r+1/n)^d r^{-d-1}dr
\geq \int_{1/n}^1 r^{-1}dr=\log n. \nonumber
\end{eqnarray}
Next, let $k=(k_1,\ldots,k_d)$ be such that $\bar{X}_2=\bar{X}_1+k/n$, and 
define $k_{\max} :=\max\{|k_1|,\ldots,|k_d|\}.$
We get that for $k_{\max} \geq 2,$
\begin{eqnarray} \label{eqn:muR1R2}
\lefteqn{\mu(R_1^n \cap R_2^n)
=\int_{(k_{\max}-1)/n}^1\int_{\BR^d}I(X(x,r)\cap \bar{X}_1\neq \emptyset,
X(x,r)\cap \bar{X}_2\neq \emptyset)dx \nu(dr) }\\
& & \leq \int_{(k_{\max}-1)/n}^1\int_{\BR^d}
I(X(x,r)\cap \bar{X}_1\neq \emptyset)dx r^{-d-1}dr
=\int_{(k_{\max}-1)/n}^1(r+1/n)^d r^{-d-1}dr \nonumber \\
& & \leq -\log((k_{\max}-1)/n)+2^d , \nonumber
\end{eqnarray}
where the last inequality follows by using the calculations in
\eqref{eqn:ref1} combined with the upper bound of \eqref{eqn:ref2}. 
Therefore, if $\bar{X}_1\neq \bar{X}_2$ and $k_{\max} \geq 2,$
we have that by using \eqref{eqn:muR1} and \eqref{eqn:muR1R2}, 
\begin{eqnarray} \label{eqn:PX1X2}
\lefteqn{\BP(\bar{X}_1,\bar{X}_2\in m_n)
=\exp\left(-\lambda \mu(R_1^n\cup R_2^n)\right)
=\exp(-2\lambda \mu(R_1^n)+\lambda \mu(R_1^n \cap R_2^n))}\\
& & \leq e^{-2\lambda \log n+\lambda 2^d-\lambda\log((k_{\max}-1)/n) }
= n^{-2\lambda}e^{\lambda 2^d}((k_{\max}-1)/n)^{-\lambda}
=e^{\lambda 2^d}(n(k_{\max}-1))^{-\lambda}. \nonumber
\end{eqnarray}
If however $k_{\max}\leq 1,$ then we simply use that 
\[
\BP(\bar{X}_1,\bar{X}_2\in m_n)\leq \BP(\bar{X}_1\in m_n)
=\BP([0,1/n]^d\in m_n).
\]
Thus, by \eqref{eqn:Pmn} and \eqref{eqn:PX1X2},
\begin{eqnarray*}
\lefteqn{\BE[m_n^2]=\sum_{\bar{X}_1\in \bar{\CX}_n}
\sum_{\bar{X}_2\in \bar{\CX}_n}\BP(\bar{X}_1,\bar{X}_2\in m_n)}\\
& & \leq n^d \left(3^d \BP([0,1/n]^d\in m_n)
+\sum_{k_{\max}=2}^n 2d k_{\max}^{d-1}e^{\lambda 2^d}(n(k_{\max}-1))^{-\lambda}\right) \\
& & \leq 3^d n^{d-\lambda}+2^d de^{\lambda 2^d}n^{d-\lambda}
\sum_{k_{\max}=2}^n (k_{\max}-1)^{d-1-\lambda}.
\end{eqnarray*}
Here, the first inequality uses that there are $n^d$ possible 
choices of $\bar{X}_1,$ and given the choice of $\bar{X}_1$
there are at most $3^d$ choices of $\bar{X}_2$ that are either the
same as, or immediate neighbours to, $\bar{X}_1.$ The remaining 
boxes $\bar{X}_2$ have $k_{max}\geq 2.$

We see that if $\lambda<d,$ then there exists a $C=C(\lambda)>0$ such that
$\BE[m_n^2]\leq C n^{2(d-\lambda)}$. Using \eqref{eqn:fmoment}
we conclude that
\[
\BP(m_n>0)\geq \frac{\BE[m_n]^2}{\BE[m_n^2]}
\geq \frac{\left(e^{-\lambda2^d}n^{d-\lambda}\right)^2}{C n^{2(d-\lambda)}}\geq c,
\]
as desired.
\fbox{}\\

Our next lemma gives a useful consequence of $\CC(\lambda)$ surviving,
but first we need some more notation. Let $D_n=D_n(\CC_n)$ be a minimal 
collection of boxes in $\bar{\CX}_n$ such that 
\[
\CC_n\cap [0,1]^d \subset \bigcup_{\bar{X}\in D_n} \bar{X}.
\]
Note that $D_n$ is not necessarily unique, as a point $x\in \CC_n$
sitting on the boundary between two boxes $\bar{X}_1$ and 
$\bar{X}_2$ can be covered by either one of them.
If there is more than one way of choosing such a set $D_n,$ 
we pick one according to some predetermined rule.
Let $L_n=|\{\bar{X}\in \bar{\CX}_n: \bar{X}\in D_n\}|$.

\begin{lemma} \label{lemma:Dninfty}
Let $\CC(\lambda)$ be either $\CC^{ball}$ or $\CC^{box}$.
For any $\lambda>0$ we have that
\[
\BP(\{\CC(\lambda)\cap[0,1]^d\neq \emptyset\} \setminus
\{\lim_{n \to \infty} L_{2^n}=\infty\})=0.
\]
\end{lemma}
{\bf Remarks:} 
The reason for proving Lemma 
\ref{lemma:Dninfty} along a subsequence $(2^n)_{n\geq 1}$, is that this 
will avoid unnecessary technical details. It is also all that we need in 
order to prove Theorem \ref{thm:main}.

Observe that if $\lim_{n \to \infty} L_{2^n}=\infty,$ then
$\CC_n\cap[0,1]^d\neq \emptyset$ for every $n\geq 1.$ As above, it follows that
also $\CC \cap[0,1]^d \neq \emptyset$.  

\bigskip

\noindent
{\bf Proof.}
Let
\[
E_{2^n}:=\bigcup_{\bar{X} \in D_{2^n}} \bar{X} ,
\]
and observe that by definition of $D_{2^n},$ we have that
\[
(\CC_{2^n}\cap[0,1]^d)\setminus E_{2^n}=\emptyset.
\]
We have that for some $\alpha=\alpha(\lambda)>0,$
\[
\BP(\CC_2\cap [0,1]^d=\emptyset)=\alpha.
\]
By using the FKG inequality for Poisson processes together
with the semi scale invariance of the models, we conclude that
\begin{eqnarray*}
\lefteqn{\BP(\CC_{2^{n+1}}\cap [0,1]^d=\emptyset | D_{2^n})}\\
& & \geq \BP(\CC_{2^{n+1}}^{2^n}\cap E_{2^n}=\emptyset | D_{2^n})
\geq \prod_{\bar{X}\in D_{2^n}}
\BP(\CC_{2^{n+1}}^{2^n}\cap \bar{X}=\emptyset)
= \alpha^{L_{2^n}}>0.
\end{eqnarray*}
Therefore, if there exists $L<\infty$ such that
$L_{2^n}\leq L$ for infinitely many $n,$
we can use L\'evy's Borel-Cantelli lemma, to conclude that almost surely
$\CC\cap [0,1]^d=\emptyset$.
\fbox{}\\

We can now prove our main result.

\noindent
{\bf Proof of Theorem \ref{thm:main}.}
The fact that $\lambda_e^{ball}=d/v_d$ is an immediate consequence
of Theorem \ref{thm:BE} as explained in Section \ref{sec:defmodels}.
Furthermore, Proposition \ref{prop:crit}
shows that $\lambda_e^{box}\geq d.$ Therefore, it remains to prove that
$\lambda^{box}_e\leq d$ and that both the ball and the box
models are in the empty phase
at their respective critical points.

Obviously, if $\bar{X}\in D_n,$ then $\bar{X}$ cannot be covered by a single
set in the Poisson process $\Phi_n$. Therefore, 
\begin{equation} \label{eqn:DnNn}
L_n\leq |M_n|.
\end{equation}

We proceed by bounding $\BE[|M_n|]$ in the case $\CC=\CC^{box}.$
Similar to the proof of Proposition \ref{prop:crit}
we have that for any $\bar{X}\in \bar{\CX}_n,$
\begin{eqnarray*}
\lefteqn{\BP(\bar{X}\in M_n)
=\exp\left(-\lambda\int_{1/n}^1 \int_{\BR^d} I(x\in(-r/2+1/n,r/2)^d)dx \nu(dr)\right)}\\
& & =\exp\left(-\lambda\int_{1/n}^1 r^{-d-1}
\sum_{k=0}^d {d \choose k} r^k (-n)^{k-d} dr\right)\\
& & =\exp\left(-\lambda\log n-
\lambda\sum_{k=0}^{d-1} {d \choose k} (-n)^{k-d}
\left(\frac{n^{d-k}-1}{d-k}\right)\right).\\
\end{eqnarray*}
Furthermore, since
\[
\sum_{k=0}^{d-1} {d \choose k} (-n)^{k-d}
\left(\frac{n^{d-k}-1}{d-k}\right)
\geq -\sum_{k=0}^{d} {d \choose k}=-2^d,
\]
we conclude that
\begin{equation}\label{eqn:emn}
\BE[|M_n|]=n^d\BP([0,1/n]^d\in M_n)\leq n^{d-\lambda}e^{\lambda 2^d}.
\end{equation}
By \eqref{eqn:DnNn},\eqref{eqn:emn} together with Lemma \ref{lemma:Dninfty}, 
we have that if $\BP(\CC^{box}\cap[0,1]^d\neq \emptyset)>0,$ then
\[
\lim_{n \to \infty}e^{\lambda 2^d}\left(2^n\right)^{d-\lambda}
\geq \lim_{n \to \infty}\BE[|M_{2^n}|]
\geq \lim_{n \to \infty}\BE[L_{2^n}]=\infty,
\]
and so we conclude that we must have $\lambda<d.$ This proves that
$\lambda_e^{box} \leq d$ and that for $\lambda=d$
\[
\BP(\CC^{box}(\lambda)\cap[0,1]^d\neq \emptyset)=0.
\]

We now turn to the case of $\CC^{ball}.$ First we observe that
for any $K\in \Phi,$
$[0,1/n]^d\subset K$ iff the closed ball
$\bar{B}(\frac{1}{2n}(1,\ldots,1),\sqrt{d}/(2n))\subset K,$
simply because of the fact that the
sets $K\in \Phi$ are balls.
We then get that for some constant $C=C(\lambda)<\infty,$ and $n>\sqrt{d}/2$
\begin{eqnarray*}
\lefteqn{\BP(\bar{X}\in M_n) 
= \BP( \not \exists K\in \Phi: \bar{B}(o,\sqrt{d}/(2n))\subset K)}\\
& & =\exp\left(-\lambda \mu(\{(x,r):\bar{B}(o,\sqrt{d}/(2n))\subset B(x,r)\})\right)\\
& & =\exp\left(-\lambda \int_{\sqrt{d}/(2n)}^1 \int_{\BR^d}
I(|x|\leq r-\sqrt{d}/(2n))dx \nu(dr)\right)\\
& & =\exp\left(-\lambda \int_{\sqrt{d}/(2n)}^1 v_d (r-\sqrt{d}/(2n))^d
r^{-d-1}dr\right)\\
& & =\exp\left(-\lambda v_d \int_{1/n}^{2/\sqrt{d}} (s-1/n)^d
s^{-d-1}ds\right)
\leq Cn^{-\lambda v_d},
\end{eqnarray*}
where the last inequality follows as above.
As for the box model, we obtain that for $\lambda_e=d/v_d$
\[
\BP(\CC^{ball}(\lambda_e)\cap [0,1]^d \neq \emptyset)=0.
\]
\fbox{}\\
\noindent
{\bf Acknowledgements} The authors would like to thank the anonymous 
referees for providing helpful suggestions improving the readability of the 
paper.


\begin{thebibliography}{99}

\bibitem{BC} Broman E. and Camia F.
Universal behavior of connectivity properties in fractal percolation models.
{\em Electron. J. Probab.} {\bf 15}, (2010), 1394--1414.

\bibitem{BE} Bierm\'{e} H. and Estrade A. Covering the whole space with
Poisson random balls. {\em ALEA Lat. Am. J. Probab. Math. Stat.}
{\bf 9}, (2012), 213--229.

\bibitem{DM} Dekking, F. M. and Meester, R. W. J.
On the structure of Mandelbrot's percolation process and other random Cantor sets.
{\em J. Statist. Phys.} {\bf 58}, (1990), no. 5-6, 1109--1126.

\bibitem{LW} G.F. Lawler and W. Werner, The Brownian loop soup.
{\em Probab. Theory Relat. Fields}, {\bf 128}, (2004), 565--588.

\bibitem{Shepp} Shepp L.A. Covering the line with random intervals
{\em Z. Wahrscheinlichkeitstheorie und Verw. Gebiete}, {\bf 23}, (1972), 163--170.

\end{thebibliography}
\end{document}